\documentclass{elsart}

\usepackage{amssymb}

\begin{document}
\newcommand{\alpt}{\alpha_1, \ldots , \alpha_t}
\newcommand{\R}{{\mathbb{R}}} 
\renewcommand{\P}{\mathfrak{P}}         
\newcommand{\p}{\mathfrak{P}}         
\newcommand{\N}{\mathbb{N}}             
\newcommand{\Z}{\mathbb{Z}}             
\newcommand{\T}{\mathbb{T}}             
\newcommand{\eps}{\varepsilon}
\newcommand{\Q}{\mathbb{Q}}
\renewcommand{\Cup}{\bigcup}
\renewcommand{\Cap}{\bigcap}
\renewcommand{\subset}{\subseteq}
\renewcommand{\supset}{\supseteq}
\newcommand{\mal}{\!\cdot\!}  
\newcommand{\F}{\mathfrak{F}}  
\newcommand{\app}{^{\smallfrown}}
\newcommand{\Vee}{\bigvee}
\newcommand{\Wedge}{\bigwedge}
\renewcommand{\tilde}{\widetilde}
\newcommand{\cddot}{}

\begin{frontmatter}
\title{Strong characterizing sequences of countable groups}
\author{Mathias Beiglb\"ock\thanksref{FWF}}\ead{mbeigl@osiris.tuwien.ac.at}
\address{Institute for Analysis and Technical Mathematics,
Vienna University of technology, 1040 Vienna, Austria}
\thanks[FWF]{Supported by the FWF-project 8312}
%
%
\begin{abstract}
 Andr\'as Bir\'o and Vera S\'os 
 prove  that for any
subgroup $G$ of $\T$ generated freely by finitely many generators
there is a sequence $A\subset \N$ such that for all $\beta \in \T$
we have ($\|.\|$ denotes the distance to the nearest integer)
$$
\beta\in G  \Rightarrow  \sum_{n\in A} \| n \beta\| <
\infty,\quad \quad \quad \beta\notin G  \Rightarrow  \limsup_{n\in
A, n \rightarrow \infty } \|n \beta\| > 0. $$ We extend this
result to arbitrary countable subgroups of $\T$. We also show
that not only the sum of norms but the sum of arbitrary small
powers of these norms can be kept small. Our proof combines ideas
from the above article with new methods, involving a filter
characterization of subgroups of $\T$.
\end{abstract}
\end{frontmatter}
%
%
\section{Introduction}
 We study certain subgroups of $\T=\R/\Z$ and methods to describe
them by sequences of positive integers. By $\|.\|$ we denote the
distance to the nearest integer. It is easily seen that for any
sequence $A\subset \N$ the set $\{ \beta\in\T : \lim_{n\in A, n
\rightarrow \infty} \| n \beta \| =0\}$ is a subgroup of $\T$. It
seems  natural to ask which subgroups arise in this way. In
\cite{bds} A. Bir\'o, J.-M. Deshouillers and V. T. S\'os show
that for any countable group $ G < \T$ there is some $A\subset
\N$ that characterizes $G$ in the above sense.
\\ Another way to connect subsets of $\N$ and $\T$ is to consider
the set $\{ \beta\in\T: \sum_{n\in A} \| n\beta \| <\infty\}$
which again is a subgroup of $\T$. Following a question of P.
Liardet  A. Bir\'o and V. T. S\'os show in \cite{bs} that if
$1,\alpt \in\T$ are linearly independent over the rationals there
is a sequence $A\subset \N$, that characterizes $\langle \alpt
\rangle $ simultaneously in both ways. Such a sequence is called
a `strong characterizing sequence' of $\langle \alpt \rangle$.
Our aim  is to find strong characterizing sequences for arbitrary
countable subgroups of $\T$. The main result is:
\begin{thm}\label{hauptsatz}
  Let $G= \{\alpha_t:t\in \N\} $ be a subgroup of $\T$. Then there
  exists a sequence $A\subset \N$, such that for all $\beta \in
  \T$
 $$
\beta\in G  \Rightarrow \forall r>0 \ \sum_{n\in A} \| n
\beta\|^r < \infty, \quad\quad\quad \beta\notin G  \Rightarrow
\limsup_{n\in A, n \rightarrow \infty} \|n \beta\| \geq 1/6.
$$
\end{thm}
%
%
\section{Connecting two methods}
In our proof we use the following reformulation of Theorem 1 in
\cite{w}:
\begin{prop}\label{reinhardsatz}
Let $G$ be an arbitrary subgroup of $\T$. Then there is a
 filter $\F$ on $\N$
  that characterizes  $G$ in the  sense that that for all
  $\beta\in\T$
$$ \beta \in G \quad \Longleftrightarrow \quad \F - \lim_{
n} \| n \beta \| =0.$$
\end{prop}
Here `$\F - \lim_{ n} \| n \beta \| =0$' means that for all
$\eps>0$ one has $\{n\in \N : \| n\beta\| \leq \eps \}\in \F$.
The filter-convergence defined in this way is more general than
ordinary convergence: For a sequence $A\subset \N$ let $\F(A)$ be
the filter consisting of all sets containing $\{ k\in A: k\geq
n\}$ for some $n\in \N$. Then we have for all $\beta \in \T$
$$ \lim_{n\in A, n\rightarrow \infty} \| n\beta\| =0  \quad \Longleftrightarrow \quad \F(A) - \lim_{n} \| n\beta\| =0.$$
The following notation will be useful: Given $\alpt\in\T$, $\eps
> 0$ and $N\in \N$ the corresponding infinite respectively finite  Bohr sets are defined by
\begin{eqnarray*}
H_\eps (\alpt)&:= & \{ n\in \N : \|n \alpha_1\| , \ldots ,\|n \alpha_t\| \leq
\eps\},\\
H_{N,\eps} (\alpt)&:= & \{ n\leq N : \|n \alpha_1\| , \ldots ,\|n \alpha_t\| \leq
\eps\}.
\end{eqnarray*}
 Using the finite intersection property of filters, one sees
that $\F - \lim_{ n} \| n \beta \| =0$ for all elements $\beta$
of some given $G<\T$ implies that for all $\alpt \in G$ and
$\eps>0$ $H_\eps (\alpt)  \in \F$. For each subgroup $G<\T$ there
is a canonical (i.e. smallest) candidate for a filter that
 characterizes $G$, namely the filter $\F_G$  which consists of all
sets containing a set $H_\eps (\alpt)$ $ (\eps>0, t\in \N, \alpt\in G)$.\\
To illustrate the connections between the number theoretic
approach in \cite{bds} respectively \cite{bs} and the more
abstract point of view in \cite{w} we show that the result on the
characterization of countable subgroups by sequences of positive
integers in \cite{bds} implies Proposition \ref{reinhardsatz}:
\begin{pf}{}Let $G<\T$ be an arbitrary subgroup and let $\F_G $ be
the filter described above. By definition of $\F_G$ we have $\F_G
- \lim_{n} \| n\beta \| =0$ for each $\beta \in G$. Now assume
$\F_G - \lim_{n} \| n \beta \| =0 $ for some $\beta\in \T$. For
$k\in \N$ let $M_k:=H_{1/k}(\beta)\in \F_G$. According to the
construction of $\F_G$, there are sequences $t_1<  t_2< \ldots\
(t_k\in\N)$, $ (\alpha_t)_{t\in \N}\ (\alpha_t\in G)$ and
$\eps_1> \eps_2
> \ldots\ (\eps_k >0)$ such that  $M_k \supset
H_{\eps_k}(\alpha_1,
\ldots , \alpha_{t_k})$ for all $k\in \N$.\\
By the result of A. Bir\'o, J.-M. Deshouillers and V. T. S\'os
there is a sequence $A\subset \N$, such that $$\{\beta \in \T:
\lim_{n\in A, n\rightarrow \infty} \| n\beta \| = 0\} = \langle
\alpha_t: t\in \N \rangle .$$
 In particular we have $\lim_{n\in A,
n\rightarrow \infty } \| n \alpha_t \| =0$ for all $t\in \N$.
Thus for fixed $m\in \N$ we can find $n_m \in \N$  satisfying $\|
n \alpha_t\|\leq \eps_m $ for all $n\in A, n\geq n_m$ and for all
$t\leq t_m$. This implies $ \{n\in A: n> n_m\}\subset H_{\eps_m
}(\alpha_1,\ldots , \alpha_{t_m})\subset M_m$, i.e. for all $n\in
A, n\geq n_m $ we have $\|n \beta\| \leq 1/m$. Since $m\in \N$
was arbitrary this yields $ \lim_{n\in A, n\rightarrow \infty} \|
n \beta \| =0$ and, as $A$ is a characterizing sequence, $\beta
\in \langle \alpha_t :t\in \N \rangle < G$.\qed
\end{pf}
%
%
\section{Ideas of the proof}
The rest of this article focuses on the proof of Theorem
\ref{hauptsatz}. The proof splits in several lemmas. Before we
state and prove them rigorously, we want to give a short sketch of
the strategy of the proof and the informal meaning of the
individual lemmas:
\\ Lemma \ref{liminfdarstellung} shows how the countable group $G$
may be represented as the limes inferior of certain open subsets
$V_t$ of $\T$. These sets may by seen as approximations of $G$.
\\ Lemma \ref{zahlentheorielemma} shows that the behaviour of the
values $\| n\beta \|$, where $n$ runs in an appropriate finite
Bohr set, may decide whether $\beta $ lies in an approximation
$V_t$ of $G$. Part (1) of the Lemma uses Theorem
\ref{reinhardsatz}, while part (2) follows easily by a
compactness argument similar to the reasoning in  \cite{bds}.
\\ The methods developed so far are powerful enough to prove the existence of
 sequences that characterize countable groups in the sense of
 \cite{bds}. To provide a strong characterizing
 sequence we use Lemma \ref{meinlemma} to replace a Bohr set $H$
 by a
 somewhat thinner  set $S$ that contains the
 same amount of information but  allows in addition to keep
 the sum $\sum_{n\in S} \| n \alpha \|^r\ (\alpha \in G, r>0)$ under
 control. The proof of Lemma \ref {meinlemma} is based on Lemma
 \ref{boormengenstruktur}, a deep result on the structure of Bohr
 sets due to A. Bir\'o and V.
T. S\'os (\cite{bs}).

%
%
\section{Preparations}
The following technical facts will be needed later. The proof is
elementary, so we skip it.
\begin{lem}\label{abschaetztrick}Let $\alpha, \beta \in\T$ and $n\in \N$.
\begin{itemize}
\item[(1)] Assume $\|\alpha\|, \| 2 \alpha \|, \ldots , \| n
\alpha \| \leq d <1/3$. Then $\| \alpha\|\leq d/n$.
\item[(2)] Assume $\| \beta + 2^0 \alpha \|,\| \beta + 2^1 \alpha \|, \ldots , \|
\beta+2^n \alpha \| \leq d <1/6$. Then $\| \alpha\|\leq
d/2^{n-2}$.
\end{itemize}
\end{lem}
Given $\alpt\in \T$ and $M\in \N$ we define $$\langle
\alpt\rangle_M :=\{ k_1\alpha_1 + \ldots + k_t \alpha_t : |k_1|,
\ldots , |k_t| \leq M\}. $$ We further define $ \| \beta S\| :=
\sup \{ \| n\beta\|: n\in S\}$ for $\beta \in \T$ and $S\subset
\N$.
\begin{lem}\label{zahlentheorielemma}
Let $\alpt\in \T$ and $\eps >0$. \begin{itemize} \item [(1)]
 There
exists some positive integer $M$ such that
$$ \| \beta H_\eps (\alpt)\| \leq 1/6 \ \Rightarrow \ \beta \in
\langle \alpt\rangle_M.$$
\item[(2)]
 If $V\supset\langle \alpt\rangle_M$ is an open subset of
$\T$, there exists some positive integer $N$ such that
$$\| \beta H_{N,\eps} (\alpt)\| \leq 1/6 \ \Rightarrow \ \beta \in
V.$$
\end{itemize}
\end{lem}
\begin{pf}
Throughout the proof we suppress mentioning $\alpt $ while notating
Bohr   sets.
\begin{itemize}
\item[(1)]
Suppose $\beta$ satisfies $\| \beta H_\eps \| \leq 1/6$. Let $\F
$ be a filter on $\N$ that characterizes $\langle \alpt\rangle$
and let $m\in \N$ be  fixed. Of course we have $H_{\eps /m} \in
\F$. For $n\in H_{\eps/m}$ and $k\leq m$ we have $k  n\in
H_{\eps}$ and in particular $\|k\cddot n \beta \| \leq 1/6$. Since
$1/6< 1/3$ this implies $\| n\beta \| \leq \frac{1}{6m}$ by Lemma
\ref{abschaetztrick}. Thus we have $\| \beta H_{\eps/m}\| \leq
\frac{1}{6m}$ and since $m$ was arbitrary we get $\F-\lim_{n} \|
n\beta \| = 0$. $\F$ was assumed to characterize $\langle
\alpt\rangle$ thus we have $\beta \in\langle
\alpt\rangle$.\\
It remains to show that $\{\beta \in \T: \| \beta H_\eps \|\leq
1/6\}$ is finite. The torsion subgroup of $\langle \alpt\rangle$
is finite and cyclic, let its order be $q\in\N$. Then $q \langle
\alpt\rangle $ is torsion free, hence we find some $\gamma_1,
\ldots , \gamma_n\in \T$, such that $q\langle \alpt\rangle$ is
freely generated by $q \gamma_1, \ldots , q \gamma_n$. We have
$\langle \alpt\rangle=\langle \gamma_1, \ldots, \gamma_n, 1/q
\rangle$ and there are uniquely determined $k_{i j}\in \Z\ (i\leq
t,j\leq n)$ and $k_i\in \{0,\ldots , q-1\}\ (i\leq t)$, such that
$$ \alpha_i = \sum_{j=1}^n k_{ij} \gamma_{j} + k_i  / q\ (i\leq t).$$
Thus we can find some $\delta >0$, such that for all $m\in q \N $
$$ \|m \gamma_1\|, \ldots, \|m \gamma_n\| \leq \delta \ \Rightarrow
\| m \alpha_1\|, \ldots, \| m \alpha_t\| \leq \eps.$$ For each
$\beta $ satisfying $\| \beta H_\eps\|\leq 1/6$  there are
uniquely determined $k_j\in \Z\ (j\leq n)$ and $k\in \{0,\ldots,
q-1\}$ such that $ \beta = \sum_{j=1}^n k_j \gamma_j+ k/q$. If
the $k_j\ (j\leq n)$ don't vanish simultaneously, Kronecker's
theorem assures that we can find $m\in q \N$, such that
\begin{eqnarray*}
  \forall j\leq n \quad  \frac{1}{6\sum_{i=1}^n |k_i| }
\quad<& \mbox{sign} (k_j) m \gamma_j &<\quad
\frac{5}{6\sum_{i=1}^n |k_i| } \quad \mbox{mod}\ 1 \\
  \Longrightarrow \qquad \quad  \frac{1}{6}\ \quad \qquad <& m
\sum_{i=1}^n k_i   \gamma_i &<\ \qquad \quad\frac{5}{6} \quad\quad
\quad  \mbox{mod}\ 1,
  \end{eqnarray*} i.e. $\| m\beta \| >1/6.$  Thus $\frac{5}{6\sum_{i=1}^n |k_i|
  }> \delta$. This shows that there are only finitely many choices
  for the $k_i\ (i\leq n)$. Thus $\{\beta \in \T: \| \beta H_\eps \|\leq
  1/6\}$ is also finite and we can find
  some $M\in \N$, such that $\{\beta \in \T: \| \beta H_\eps \|\leq
  1/6\}\subset \langle\alpt\rangle_M$.
\item[(2)] Let $M$ be as in (1). Then $\langle \alpt \rangle_M \subset V$ implies
  $$ \emptyset = V^c \cap \{\beta \in \T:\| \beta H_\eps\|\leq 1/6\}=
  V^c \cap \Cap_{n\in H_\eps } \{\beta \in \T:\| n \beta\|\leq
  1/6\}.$$  Since $\T$ is compact and all of the above sets are
  closed, the intersection of finitely many of these sets
  must be
  empty, i.e. we can find some $N\in \N $ such that $V^c \cap \Cap_{n\in H_{N,\eps} } \{\beta \in \T:\| n \beta\|\leq
  1/6\}=\emptyset$. Obviously this $N$ is as required.
\end{itemize}\qed
\end{pf}
\begin{lem}\label{liminfdarstellung}
Let $G=\{ \alpha_t: t\in \N\}$ be a subgroup of $\T$  and let
$(M_t)_{t\in \N}$ be a sequence of positive integers. There exists
a sequence $(V_t)_{t\in \N}$ of open subsets of $T$ such that
\begin{itemize}
\item[(i)] $V_t \supset \langle \alpt\rangle_{M_t}\ (t\in \N),$
\item[(ii)] $ \Cup_{k\in \N} \Cap_{t\geq k} V_t =
\liminf_{t\rightarrow \infty} V_t =G.$
\end{itemize}
\end{lem}
\begin{pf}
We may assume that $(M_t)_{t\in\N}$ is increasing. We choose a
sequence $(\delta_t)_{t\in \N}$ of positive numbers that decreases
to $0$ and satisfies for all $t\in \N$
\begin{itemize}
 \item[(1)] $2 \delta_t < \min \{\|\alpha-\alpha' \|:\alpha, \alpha' \in
\langle \alpt \rangle_{M_t}, \alpha \neq
 \alpha'\}$,
 \item[(2)] $\delta_t+ \delta_{t+1} < \min \left\{\|\alpha-\alpha'
\|:\begin{array}{cl} &\alpha \in \langle
 \alpt
 \rangle_{M_t},\\
  &\alpha' \in \langle \alpha_1,\ldots ,\alpha_{t+1} \rangle_{M_{t+1}} \setminus \langle
  \alpt
 \rangle_{M_t}
 \end{array}\right\}.$
\end{itemize} Using this, we define
$$V_t := \{\beta\in \T :\exists \alpha\in \langle \alpt\rangle_{M_t}
 \ \| \alpha-\beta\| < \delta_t\}.$$ We obviously have
$\liminf_{t\rightarrow \infty } V_t \supset G$. To show the
reverse inclusion, assume $\beta \in \liminf_{t\rightarrow \infty
} V_t$, i.e. $\beta\in V_t$ for all $t\geq t_0$ for some
$t_0\in\N$. By definition of the $V_t$ for all $t\geq t_0$ there
is some $\gamma_t\in \langle \alpt\rangle_{M_t}$ satisfying $\|
\beta- \gamma_t \| < \delta_t$ and (1) shows that this $\gamma_t$
is uniquely determined. Further $\gamma_t \neq \gamma_{t+1}$ for
some $t\geq t_0$ would contradict (2), thus we have
$\gamma_{t_0}= \gamma_{t_0+1}= \gamma_{t_0+2} = \cdots$. In
particular this shows $ \| \beta- \gamma_{t_0} \| =\| \beta -
\gamma_t\| < \delta_t \rightarrow 0$, hence $\beta =
\gamma_{t_0}\in G$. \qed
\end{pf}

From Lemma 1 in \cite{bs} one gets:
\begin{lem}\label{boormengenstruktur}
  Let $t\in \N$. There exists some constant $C_1=C_1(t)$,
   such that for all  $ \alpha_1, \ldots , \alpha_t\in \T$, positive $\eps\leq 1/C_1$
  and positive integers $N$ there are suitable
  nonzero integers $n_1, \ldots , n_R $ and
  positive integers $K_1, \ldots , K_R$, $R\leq C_1$ satisfying
  \begin{itemize}
  \item[(a)] $\sum_{i=1}^R K_i \| n_i \alpha_j\| \leq  C_1 \cddot \eps \quad ( 1 \leq j\leq t) $
  \item [(b)] $ \sum_{i=1}^R K_i | n_i| \leq C_1 \cddot N$,
  \item[(c)] $ H_{N,\eps}(\alpt) \subset \left\{ \sum_{i=1}^R k_i n_i : 1\leq k_i \leq K_i\right\}$.
  \end{itemize}
\end{lem}

\begin{lem}\label{meinlemma}  Let $t\in \N$. There exists some constant
$C_2=C_2(t)$,
   such that for all  $ \alpha_1, \ldots , \alpha_t\in \T$, positive $\eps\leq 1/C_1(t)$, positive $r\leq 1$ and
   positive integers $N$ and $U$ there is a suitable nonempty finite set $S$ of integers satisfying
   \begin{itemize}
   \item[(i)] $U < \min S$,
   \item[(ii)] for all $j\leq t$ we have  $\sum_{n\in S} \| n \alpha_j\|^r \leq C_2 \cddot \frac{\eps^r}{2^r-1}, $
   \item[(iii)] for all $\beta \in \T$ we have $ \min \{1/6,\|\beta  H_{N,\eps}(\alpt) \| \} \leq \| \beta S \|   $.
   \end{itemize}
\end{lem}
\begin{pf} Let $\alpt\in \T$ and $C_1,R, K_i, n_i \ (i \leq R) $  as given by Lemma \ref{boormengenstruktur}.
Let $m > U$ be an integer satisfying $$\| m \alpha_j\|^r \leq
\frac{\eps^r}{\lg_2 (8 \cddot C_1^2\cddot N )}$$ for all $j\leq t$
and let
$$ S=  \{m + 2^l |n_i| : 2^l\leq 8\cddot K_i\cddot R\}. $$
Clearly $S$ satisfies (i).\\
For each $j\leq t$ we have
$$ \sum_{n\in S } \|n \alpha_j \|^r \leq  \mbox{card} (S) \cddot \| m \alpha_j \|^r +
\sum_{n\in S} \| (n-m) \alpha_j\|^r.$$ To find an upper bound for
the first term, we observe that $K_i \leq C_1 \cddot N $
  implies $\mbox{card} (S) \leq R \lg_2 (8\cddot  C_1\cddot N \cddot R)$. Thus
  $$\mbox{card} (S) \cddot \| m \alpha_j \|^r \leq R \lg_2 (8\cddot C_1\cddot N \cddot R) \cddot\frac{ \eps^r }{ \lg_2 (8\cddot C_1^2 \cddot N
  )}
  \leq C_1\cddot\eps^r.$$
  The second term can be estimated by \begin{eqnarray*}
  \sum_{i=1}^R \left( \sum_{l=0}^{ \lfloor \lg_2(8\cddot K_i\cddot R)\rfloor} 2^l
\| n_i \alpha_j\|\right)^r&\leq&
   \sum_{i=1}^R \frac{\left( 2^r\right)^{\lg_2 (8 K_i R)+1}-1}{2^r-1} \| n_i
   \alpha\|^r\\
      &<& \frac{16^r R^r}{2^r-1} \sum_{i=1}^R K_i^r \| n_i
      \alpha\|^r.\end{eqnarray*}
    For any $a_1, \ldots , a_R$ we have
   $ \frac{1}{R} \sum_{i=1}^R a_i^r \leq  \left(\frac{1}{R} \sum_{i=1}^R a_i \right)^r $ by Jensen's
   inequality. This yields
   $$ \sum_{n\in S} \| (n-m) \alpha_j\|^r \leq \frac{16^r R}{2^r-1}\left( \sum_{i=1}^R K_i \|  n_i \alpha_j\| \right) ^r
   \leq  \frac{16^r C_1}{2^r-1} \left( C_1 \eps \right) ^r .
     $$
    Thus $S$ will satisfy (ii) if we let $C_2 := C_1 + 16\cddot C_1^2  $ .\\
Finally let $\beta \in \T$ and $d:= \| \beta S \|$. We may assume
$d<1/6$. Thus  by Lemma \ref{abschaetztrick} for all $i\leq R$
$$ \| m \beta + 2^l |n_i| \beta \| \leq d \ (l \leq \lg_2 (8\cddot K_i \cddot R)) $$ implies $$\|n_i
\beta \| \leq  \frac{d }{2^{\lfloor lg_2 (8 \cddot K_i \cddot
R)\rfloor -2}} \leq \frac{d}{K_i\cddot R}.$$ By Lemma
\ref{boormengenstruktur} each $n \in H_{N,\eps }(\alpt)$ has a
representation $n= \sum _{i=1 }^R k_i n_i $ for some integers
$k_i, (1\leq i \leq R)$ satisfying $1\leq k_i \leq K_i$. Using
this representation we get
$$ \| n\beta \| =   \left\| \sum _{i=1 }^R k_i n_i \beta\right\|
\leq \sum _ {i= 1}^R K_i \| n_i \beta\| \leq \sum _ {i= 1}^R K_i
\frac{d }{K_i\cddot R}= d. $$ Thus $S$ satisfies (iii). \qed
\end{pf}
\section{Proof of the Theorem}
Finally we are able to give the proof of Theorem \ref{hauptsatz}.
 Let $(\eps_t)_{t\in \N}$ be a sequence of positive
numbers, satisfying  $\eps_t < 1/C_1(t)$ and $\sum_{t=1}^\infty
C_2(t) \frac{\eps_t^{1/t}}{2^{1/t}-1} <\infty$. Combining Lemma
\ref{zahlentheorielemma} and Lemma \ref{liminfdarstellung} we
find a sequence $(N_t)_{t\in \N}$ of positive integers and a
sequence $(V_t)_{t\in \N}$ of open subsets of $\T$, such that:
\begin{itemize}
\item[(1)] For all $\beta \in \T$ and for all $t\in \N$\ $\| \beta H_{N_t, \eps_t}(\alpt) \| \leq 1/6 \
\Rightarrow \ \beta\in V_t.$
\item[(2)] $ \Cup_{k\in \N} \Cap_{t\geq k} V_t  =G.$
\end{itemize}
Using Lemma \ref{meinlemma} we find some sequence $(S_t)_{t\in
\N}$ of subsets of $\N$ such that for all $t\in \T$
 \begin{itemize}
   \item[(i)] $\max S_t < \min S_{t+1}$,
   \item[(ii)] $ \sum_{n\in S_t} \| n \alpha_j\|^{1/t} \leq C_2 (t) \cddot \frac{\eps_t^{1/t}}{2^{1/t}-1} \ (  j\leq t)$,
   \item[(iii)] for all $\beta \in \T$ $ \min \{1/6, \| \beta  H_{N_t,\eps_t}(\alpt) \|\} \leq \| \beta S_t \|
   $.
   \end{itemize}
By defining $A:= \Cup_{t\in \N} S_t$ we will in fact get a strong
characterizing sequence of $G$ as stated in Theorem
\ref{hauptsatz}:
\\ Assume $\beta \in G$ and $r>0$. Then $\beta =\alpha_{t_0} $ for some
$t_0\in \N$. If we let $m> \max\{t_0,1/r\}$, we have
$$ \sum_{n\in A, n\geq \min S_{m}} \| n\beta \|^r \leq \sum_{t\geq m } \sum_{n\in S_t} \| n
\alpha_{t_0}\|^{1/t}
 \leq \sum_{t\geq t_0 } C_2(t) \frac{\eps_t^{1/t}}{2^{1/t}-1} <\infty.$$
 \\ Finally, assume $\beta \notin G$. There exists a
 sequence $t_1< t_2< \ldots $ of positive integers such that $\beta
 \notin V_{t_k}\ (k\in \N)$. So for each $k\in \N$ we have
 $\|\beta H_{\eps_{t_k},N_{t_k}} (\alpha_1,\ldots,\alpha_{t_k})
 \|>1/6$ and thus can find some $n_k\in S_{t_k}$ satisfying
 $\|\beta n_k \| \geq 1/6$. This shows $\limsup_{n\in A, n\rightarrow \infty} \|n\beta
 \| \geq 1/6.$ \qed

\begin{ack} The author would like to thank Gabriel Maresch and
Reinhard Winkler for reading the manuscript and making valuable
comments. Further the author is grateful to Andr\'as Bir\'o. His
remarks led to a significant strengthening of the main theorem.
\end{ack}

\end{document}